\documentclass[12pt]{article}
\usepackage{amssymb}

\usepackage[latin1]{inputenc}
\usepackage{longtable}
\usepackage{color}

\usepackage{url}

\begin{document}
\begin{center}
\Large
Solution for the problem of the game \emph{heads or tails} 
\end{center}

\begin{center}
Alberto Costa \\
LIX, \'Ecole Polytechnique\\
91128 Palaiseau, France\\
\verb+costa@lix.polytechnique.fr+
\end{center}
\textbf{Abstract} 
\\
In this paper, we describe the solution for a problem dealing with definite properties of binary sequences. This problem, proposed by Xavier Grandsart in the form of a mathematical contest \cite{1}, has been solved also by Maher Younan, Ph.D. student of Theoretical Physics at the University of Geneva, and Pierre Deligne, professor at Princeton and Field Medals, using different approaches with respect to the one presented in this work. All the proofs can be found in \cite{1} (in French).

\section{Introduction}
Consider a repetition of $n$ throwings in the game heads or tails: 0 if heads, 1 if tail. This can be represented as a binary sequence $D=  d_0,\,d_1,\,
\dots,\, d_{n-1}$  of length $n$.

We consider the 4 pairs made by inverse sequences of order 4 with periodical dimensions: $0100-0010, 1101-1011, 0011-1100, 1010-0101$. We compute the difference between the number of occurrences of the two terms in each pair for our binary sequence $D$. Only the inverse sequences made of consecutive terms are considered. Furthermore, the binary sequence is circular. Hence we consider only the inverse sequences of order 4 having this form:  $d_{{i}_{|n}},\,d_{{(i+1)}_{|n}},\,d_{{(i+2)}_{|n}},\,d_{{(i+3)}{|n}},\, i \in \{0,\dots,n-1\}$ (where $i_{|n}$ means $i$ mod $n$).

Prove that, for each random binary sequence $D$, the difference between the number of occurrences of the two terms in each pair of inverse sequences of order 4 is preserved. For instance, why if the number of occurrences of the term 0011 is equal to the number of occurrences of the term 1100 minus $t$, the number of occurrences of the term 1101 is equal to the number of occurrences of the term 1011 minus $t$ ?

\section{The solution}
We consider a binary random sequence $D=  d_0,\,d_1,\,
\dots,\, d_{n-1}$ with length $n$. We have to prove that, for every sequence $D$, if we consider the four pairs of inverse sequences of order 4 with periodical dimension ($0100-0010, 1101-1011, 0011-1100, 1010-0101$), the difference between the number of occurrences of the two terms in each pair is the same.
To do this, we use the mathematical induction on the size of the sequence $D$.
\\

\textsc{Theorem 1.}
 \label{po01} 
  For every $D=d_0,\,d_1,\,\dots,\,d_{n-1},\,
  d_i\in\{0,\,1\},\,0\leq i<n$, the difference between the number of occurrences of the two terms
  $0100-0010,\,1101-1011,\,1010-0101,\,0011-1100$ is always the same.
\\

 \textit{Proof}. We consider as basis of induction the sequence of size 4, then we will prove
  that for every $n>4$ the hypothesis is fulfilled.

\vspace{0.5cm}

\noindent \textbf{Basis step.}\\
There are 16 possible binary sequences $D$ of size 4. The
table \ref{tab:4kDs} shows the number of inverse sequences and the difference
between the number of occurences of the 2 terms. As can be seen, this difference is equal to 0 for each sequence.

\begin{table}[ht]
	\centering
		\begin{tabular}{|c|c|c|c|c|c|c|c|c|c|}
		\hline
		Sequence & 0100 & 0010 & {\color{red}1101} & {\color{red}1011}& {\color{green} 1010} & {\color{green}0101} & {\color{blue} 0011} & {\color{blue} 1100} & Difference \\
		\hline
		0000 & 0 & 0 & {\color{red} 0} & {\color{red} 0} & {\color{green} 0} & {\color{green} 0} &{\color{blue} 0} &{\color{blue} 0} & 0\\ \hline
		0001 & 1 & 1 & {\color{red} 0} & {\color{red} 0} & {\color{green} 0} & {\color{green} 0} &{\color{blue} 0} &{\color{blue} 0}&0\\ \hline
		0010 & 1 & 1 & {\color{red} 0} & {\color{red} 0} & {\color{green} 0} & {\color{green} 0} &{\color{blue} 0} &{\color{blue} 0}&0\\ \hline
		0011 & 0 & 0 & {\color{red} 0} & {\color{red} 0} & {\color{green} 0} & {\color{green} 0} &{\color{blue} 1} &{\color{blue} 1}&0\\ \hline
		0100 & 1 & 1 & {\color{red} 0} & {\color{red} 0} & {\color{green} 0} & {\color{green} 0} &{\color{blue} 0} &{\color{blue} 0}&0\\\hline
		0101 & 0 & 0 & {\color{red} 0} & {\color{red} 0} & {\color{green} 2} & {\color{green} 2} &{\color{blue} 0} &{\color{blue} 0}&0\\\hline
		0110 & 0 & 0 & {\color{red} 0} & {\color{red} 0} & {\color{green} 0} & {\color{green} 0} &{\color{blue} 1} &{\color{blue} 1}&0\\\hline
		0111 & 0 & 0 & {\color{red} 1} & {\color{red} 1} & {\color{green} 0} & {\color{green} 0} &{\color{blue} 0} &{\color{blue} 0}&0\\\hline
		1000 & 1 & 1 & {\color{red} 0} & {\color{red} 0} & {\color{green} 0} & {\color{green} 0} &{\color{blue} 0} &{\color{blue} 0}&0\\\hline
		1001 & 0 & 0 & {\color{red} 0} & {\color{red} 0} & {\color{green} 0} & {\color{green} 0} &{\color{blue} 0} &{\color{blue} 0}&0\\ \hline
		1010 & 0 & 0 & {\color{red} 0} & {\color{red} 0} & {\color{green} 2} & {\color{green} 2} &{\color{blue} 0} &{\color{blue} 0}&0\\ \hline
		1011 & 0 & 0 & {\color{red} 1} & {\color{red} 1} & {\color{green} 0} & {\color{green} 0} &{\color{blue} 0} &{\color{blue} 0}&0\\ \hline
		1100 & 0 & 0 & {\color{red} 0} & {\color{red} 0} & {\color{green} 0} & {\color{green} 0} &{\color{blue} 1} &{\color{blue} 1}&0\\\hline
		1101 & 0 & 0 & {\color{red} 1} & {\color{red} 1} & {\color{green} 0} & {\color{green} 0} &{\color{blue} 0} &{\color{blue} 0}&0\\ \hline
		1110 & 0 & 0 & {\color{red} 1} & {\color{red} 1} & {\color{green} 0} & {\color{green} 0} &{\color{blue} 0} &{\color{blue} 0}&0\\ \hline
		1111 & 0 & 0 & {\color{red} 0} & {\color{red} 0} & {\color{green} 0} & {\color{green} 0} &{\color{blue} 0} &{\color{blue} 0}&0\\ \hline

	\end{tabular}
	\caption{Sequences of size 4}
	\label{tab:4kDs}
\end{table}

\vspace{0.5cm}

\noindent \textbf{Inductive step.}\\
Suppose  that the theorem is valid for a sequence
$D=d_0,\,d_1,\,\dots,\,d_{n-1}$. We will show that this is also true
for the sequence $D'=d_0,\,d_1,\,\dots,\,d_{n-1},d_n$. The difference
in the inverse sequences of order 4 with periodical dimension between $D$ and
$D'$ is that, by moving from $D$ to $D'$, the following sequences are lost:
\begin{itemize}
\item $d_{n-3},\,d_{n-2},\,d_{n-1},\,d_0$
\item $d_{n-2},\,d_{n-1},\,d_{0},\,d_1$
\item $d_{n-1},\,d_{0},\,d_{1},\,d_2$
\end{itemize}

However, there are 4 new sequences:
\begin{itemize}
\item $d_{n-3},\,d_{n-2},\,d_{n-1},\,d_n$
\item $d_{n-2},\,d_{n-1},\,d_{n},\,d_0$
\item $d_{n-1},\,d_{n},\,d_{0},\,d_1$
\item $d_{n},\,d_{0},\,d_{1},\,d_2$.
\end{itemize}

Hence, the digits involved in that change are $d_{n-3}-\,d_{n-2}-\,d_{n-1}-\,d_0-\,d_1-\,d_2$ and $d_n$. Since
$d_i$ is  0 or 1, and  there are 7 digits involved, there are $2^7$ possible cases to be considered: we will show that, for every one, the
theorem 1 is valid.

Let $P$ be the sequence $d_{n-3},\,d_{n-2},\,d_{n-1},\,d_0,\,d_1,\,d_2$;
The first step to do is to compute the number of sequences
lost for each case. Table \ref{tab:spd} shows this.

	\begin{longtable}{|c|c|c|c|c|c|c|c|c|}			
\hline
	$P$	& 0100 & 0010 & {\color{red}1101} & {\color{red}1011}& {\color{green} 1010} & {\color{green}0101} & {\color{blue} 0011} & {\color{blue} 1100} \\
		\hline
		000000 & 0 & 0 & {\color{red} 0} & {\color{red} 0} & {\color{green} 0} & {\color{green} 0} &{\color{blue} 0} &{\color{blue} 0} \\ \hline
	000001	& 0 & 0 & {\color{red} 0} & {\color{red} 0} & {\color{green} 0} & {\color{green} 0} &{\color{blue} 0} &{\color{blue} 0} \\ \hline
000010& 0 & 1 & {\color{red} 0} & {\color{red} 0} & {\color{green} 0} & {\color{green} 0} &{\color{blue} 0} &{\color{blue} 0} \\ \hline
000011& 0 & 0 & {\color{red} 0} & {\color{red} 0} & {\color{green} 0} & {\color{green} 0} &{\color{blue} 1} &{\color{blue} 0} \\ \hline
000100& 1 & 1 & {\color{red} 0} & {\color{red} 0} & {\color{green} 0} & {\color{green} 0} &{\color{blue} 0} &{\color{blue} 0} \\ \hline
000101& 0 & 1 & {\color{red} 0} & {\color{red} 0} & {\color{green} 0} & {\color{green} 1} &{\color{blue} 0} &{\color{blue} 0} \\ \hline
000110& 0 & 0 & {\color{red} 0} & {\color{red} 0} & {\color{green} 0} & {\color{green} 0} &{\color{blue} 1} &{\color{blue} 0} \\ \hline
000111& 0 & 0 & {\color{red} 0} & {\color{red} 0} & {\color{green} 0} & {\color{green} 0} &{\color{blue} 1} &{\color{blue} 0} \\ \hline
001000& 1 & 1 & {\color{red} 0} & {\color{red} 0} & {\color{green} 0} & {\color{green} 0} &{\color{blue} 0} &{\color{blue} 0} \\ \hline
001001& 1 & 1 & {\color{red} 0} & {\color{red} 0} & {\color{green} 0} & {\color{green} 0} &{\color{blue} 0} &{\color{blue} 0} \\ \hline
001010& 0 & 1 & {\color{red} 0} & {\color{red} 0} & {\color{green} 1} & {\color{green} 1} &{\color{blue} 0} &{\color{blue} 0} \\ \hline
001011& 0 & 1 & {\color{red} 0} & {\color{red} 1} & {\color{green} 0} & {\color{green} 1} &{\color{blue} 0} &{\color{blue} 0} \\ \hline
001100& 0 & 0 & {\color{red} 0} & {\color{red} 0} & {\color{green} 0} & {\color{green} 0} &{\color{blue} 1} &{\color{blue} 1} \\ \hline
001101& 0 & 0 & {\color{red} 1} & {\color{red} 0} & {\color{green} 0} & {\color{green} 0} &{\color{blue} 1} &{\color{blue} 0} \\ \hline
001110& 0 & 0 & {\color{red} 0} & {\color{red} 0} & {\color{green} 0} & {\color{green} 0} &{\color{blue} 1} &{\color{blue} 0} \\ \hline
001111& 0 & 0 & {\color{red} 0} & {\color{red} 0} & {\color{green} 0} & {\color{green} 0} &{\color{blue} 1} &{\color{blue} 0} \\ \hline
010000& 1 & 0 & {\color{red} 0} & {\color{red} 0} & {\color{green} 0} & {\color{green} 0} &{\color{blue} 0} &{\color{blue} 0} \\ \hline
010001& 1 & 0 & {\color{red} 0} & {\color{red} 0} & {\color{green} 0} & {\color{green} 0} &{\color{blue} 0} &{\color{blue} 0} \\ \hline
010010& 1 & 1 & {\color{red} 0} & {\color{red} 0} & {\color{green} 0} & {\color{green} 0} &{\color{blue} 0} &{\color{blue} 0} \\ \hline
010011& 1 & 0 & {\color{red} 0} & {\color{red} 0} & {\color{green} 0} & {\color{green} 0} &{\color{blue} 1} &{\color{blue} 0} \\ \hline
010100& 1 & 0 & {\color{red} 0} & {\color{red} 0} & {\color{green} 1} & {\color{green} 1} &{\color{blue} 0} &{\color{blue} 0} \\ \hline
010101& 0 & 0 & {\color{red} 0} & {\color{red} 0} & {\color{green} 1} & {\color{green} 2} &{\color{blue} 0} &{\color{blue} 0} \\ \hline
010110& 0 & 0 & {\color{red} 0} & {\color{red} 1} & {\color{green} 0} & {\color{green} 1} &{\color{blue} 0} &{\color{blue} 0} \\ \hline
010111& 0 & 0 & {\color{red} 0} & {\color{red} 1} & {\color{green} 0} & {\color{green} 1} &{\color{blue} 0} &{\color{blue} 0} \\ \hline
011000& 0 & 0 & {\color{red} 0} & {\color{red} 0} & {\color{green} 0} & {\color{green} 0} &{\color{blue} 0} &{\color{blue} 1} \\ \hline
011001& 0 & 0 & {\color{red} 0} & {\color{red} 0} & {\color{green} 0} & {\color{green} 0} &{\color{blue} 0} &{\color{blue} 1} \\ \hline
011010& 0 & 0 & {\color{red} 1} & {\color{red} 0} & {\color{green} 1} & {\color{green} 0} &{\color{blue} 0} &{\color{blue} 0} \\ \hline
011011& 0 & 0 & {\color{red} 1} & {\color{red} 1} & {\color{green} 0} & {\color{green} 0} &{\color{blue} 0} &{\color{blue} 0} \\ \hline
011100& 0 & 0 & {\color{red} 0} & {\color{red} 0} & {\color{green} 0} & {\color{green} 0} &{\color{blue} 0} &{\color{blue} 1} \\ \hline
011101& 0 & 0 & {\color{red} 1} & {\color{red} 0} & {\color{green} 0} & {\color{green} 0} &{\color{blue} 0} &{\color{blue} 0} \\ \hline
011110& 0 & 0 & {\color{red} 0} & {\color{red} 0} & {\color{green} 0} & {\color{green} 0} &{\color{blue} 0} &{\color{blue} 0} \\ \hline
011111& 0 & 0 & {\color{red} 0} & {\color{red} 0} & {\color{green} 0} & {\color{green} 0} &{\color{blue} 0} &{\color{blue} 0} \\ \hline
100000& 0 & 0 & {\color{red} 0} & {\color{red} 0} & {\color{green} 0} & {\color{green} 0} &{\color{blue} 0} &{\color{blue} 0} \\ \hline
100001& 0 & 0 & {\color{red} 0} & {\color{red} 0} & {\color{green} 0} & {\color{green} 0} &{\color{blue} 0} &{\color{blue} 0} \\ \hline
100010& 0 & 1 & {\color{red} 0} & {\color{red} 0} & {\color{green} 0} & {\color{green} 0} &{\color{blue} 0} &{\color{blue} 0} \\ \hline
100011& 0 & 0 & {\color{red} 0} & {\color{red} 0} & {\color{green} 0} & {\color{green} 0} &{\color{blue} 1} &{\color{blue} 0} \\ \hline
100100& 1 & 1 & {\color{red} 0} & {\color{red} 0} & {\color{green} 0} & {\color{green} 0} &{\color{blue} 0} &{\color{blue} 0} \\ \hline
100101& 0 & 1 & {\color{red} 0} & {\color{red} 0} & {\color{green} 0} & {\color{green} 1} &{\color{blue} 0} &{\color{blue} 0} \\ \hline
100110& 0 & 0 & {\color{red} 0} & {\color{red} 0} & {\color{green} 0} & {\color{green} 0} &{\color{blue} 1} &{\color{blue} 0} \\ \hline
100111& 0 & 0 & {\color{red} 0} & {\color{red} 0} & {\color{green} 0} & {\color{green} 0} &{\color{blue} 1} &{\color{blue} 0} \\ \hline
101000& 1 & 0 & {\color{red} 0} & {\color{red} 0} & {\color{green} 1} & {\color{green} 0} &{\color{blue} 0} &{\color{blue} 0} \\ \hline
101001& 1 & 0 & {\color{red} 0} & {\color{red} 0} & {\color{green} 1} & {\color{green} 0} &{\color{blue} 0} &{\color{blue} 0} \\ \hline
101010& 0 & 0 & {\color{red} 0} & {\color{red} 0} & {\color{green} 2} & {\color{green} 1} &{\color{blue} 0} &{\color{blue} 0} \\ \hline
101011& 0 & 0 & {\color{red} 0} & {\color{red} 1} & {\color{green} 1} & {\color{green} 1} &{\color{blue} 0} &{\color{blue} 0} \\ \hline
101100& 0 & 0 & {\color{red} 0} & {\color{red} 1} & {\color{green} 0} & {\color{green} 0} &{\color{blue} 0} &{\color{blue} 1} \\ \hline
101101& 0 & 0 & {\color{red} 1} & {\color{red} 1} & {\color{green} 0} & {\color{green} 0} &{\color{blue} 0} &{\color{blue} 0} \\ \hline
101110& 0 & 0 & {\color{red} 0} & {\color{red} 1} & {\color{green} 0} & {\color{green} 0} &{\color{blue} 0} &{\color{blue} 0} \\ \hline
101111& 0 & 0 & {\color{red} 0} & {\color{red} 1} & {\color{green} 0} & {\color{green} 0} &{\color{blue} 0} &{\color{blue} 0} \\ \hline
110000& 0 & 0 & {\color{red} 0} & {\color{red} 0} & {\color{green} 0} & {\color{green} 0} &{\color{blue} 0} &{\color{blue} 1} \\ \hline
110001& 0 & 0 & {\color{red} 0} & {\color{red} 0} & {\color{green} 0} & {\color{green} 0} &{\color{blue} 0} &{\color{blue} 1} \\ \hline
110010& 0 & 1 & {\color{red} 0} & {\color{red} 0} & {\color{green} 0} & {\color{green} 0} &{\color{blue} 0} &{\color{blue} 1} \\ \hline
110011& 0 & 0 & {\color{red} 0} & {\color{red} 0} & {\color{green} 0} & {\color{green} 0} &{\color{blue} 1} &{\color{blue} 1} \\ \hline
110100& 1 & 0 & {\color{red} 1} & {\color{red} 0} & {\color{green} 1} & {\color{green} 0} &{\color{blue} 0} &{\color{blue} 0} \\ \hline
110101& 0 & 0 & {\color{red} 1} & {\color{red} 0} & {\color{green} 1} & {\color{green} 1} &{\color{blue} 0} &{\color{blue} 0} \\ \hline
110110& 0 & 0 & {\color{red} 1} & {\color{red} 1} & {\color{green} 0} & {\color{green} 0} &{\color{blue} 0} &{\color{blue} 0} \\ \hline
110111& 0 & 0 & {\color{red} 1} & {\color{red} 1} & {\color{green} 0} & {\color{green} 0} &{\color{blue} 0} &{\color{blue} 0} \\ \hline
111000& 0 & 0 & {\color{red} 0} & {\color{red} 0} & {\color{green} 0} & {\color{green} 0} &{\color{blue} 0} &{\color{blue} 1} \\ \hline
111001& 0 & 0 & {\color{red} 0} & {\color{red} 0} & {\color{green} 0} & {\color{green} 0} &{\color{blue} 0} &{\color{blue} 1} \\ \hline
111010& 0 & 0 & {\color{red} 1} & {\color{red} 0} & {\color{green} 1} & {\color{green} 0} &{\color{blue} 0} &{\color{blue} 0} \\ \hline
111011& 0 & 0 & {\color{red} 1} & {\color{red} 1} & {\color{green} 0} & {\color{green} 0} &{\color{blue} 0} &{\color{blue} 0} \\ \hline
111100& 0 & 0 & {\color{red} 0} & {\color{red} 0} & {\color{green} 0} & {\color{green} 0} &{\color{blue} 0} &{\color{blue} 1} \\ \hline
111101& 0 & 0 & {\color{red} 1} & {\color{red} 0} & {\color{green} 0} & {\color{green} 0} &{\color{blue} 0} &{\color{blue} 0} \\ \hline
111110& 0 & 0 & {\color{red} 0} & {\color{red} 0} & {\color{green} 0} & {\color{green} 0} &{\color{blue} 0} &{\color{blue} 0} \\ \hline
111111& 0 & 0 & {\color{red} 0} & {\color{red} 0} & {\color{green} 0} & {\color{green} 0} &{\color{blue} 0} &{\color{blue} 0} \\ \hline

	\caption{Sequences lost}
	\label{tab:spd}
\end{longtable}

After adjoining  $d_n$, there are four new sequences, as previously described. Let $d_n=0$ and $P=d_{n-3},\,d_{n-2},\,d_{n-1},\,0,\,d_0,\,d_1,\,d_2$; table
\ref{tab:sadj0} shows  the number of sequences adjoined for each
case.

	\begin{longtable}{|c|c|c|c|c|c|c|c|c|}			
\hline
	$P'$	& 0100 & 0010 & {\color{red}1101} & {\color{red}1011}& {\color{green} 1010} & {\color{green}0101} & {\color{blue} 0011} & {\color{blue} 1100} \\
		\hline
		0000000 & 0 & 0 & {\color{red} 0} & {\color{red} 0} & {\color{green} 0} & {\color{green} 0} &{\color{blue} 0} &{\color{blue} 0} \\ \hline
	0000001	& 0 & 0 & {\color{red} 0} & {\color{red} 0} & {\color{green} 0} & {\color{green} 0} &{\color{blue} 0} &{\color{blue} 0} \\ \hline
0000010& 0 & 1 & {\color{red} 0} & {\color{red} 0} & {\color{green} 0} & {\color{green} 0} &{\color{blue} 0} &{\color{blue} 0} \\ \hline
0000011& 0 & 0 & {\color{red} 0} & {\color{red} 0} & {\color{green} 0} & {\color{green} 0} &{\color{blue} 1} &{\color{blue} 0} \\ \hline
0000100& 1 & 1 & {\color{red} 0} & {\color{red} 0} & {\color{green} 0} & {\color{green} 0} &{\color{blue} 0} &{\color{blue} 0} \\ \hline
0000101& 0 & 1 & {\color{red} 0} & {\color{red} 0} & {\color{green} 0} & {\color{green} 1} &{\color{blue} 0} &{\color{blue} 0} \\ \hline
0000110& 0 & 0 & {\color{red} 0} & {\color{red} 0} & {\color{green} 0} & {\color{green} 0} &{\color{blue} 1} &{\color{blue} 0} \\ \hline
0000111& 0 & 0 & {\color{red} 0} & {\color{red} 0} & {\color{green} 0} & {\color{green} 0} &{\color{blue} 1} &{\color{blue} 0} \\ \hline
0010000& 1 & 1 & {\color{red} 0} & {\color{red} 0} & {\color{green} 0} & {\color{green} 0} &{\color{blue} 0} &{\color{blue} 0} \\ \hline
0010001& 1 & 1 & {\color{red} 0} & {\color{red} 0} & {\color{green} 0} & {\color{green} 0} &{\color{blue} 0} &{\color{blue} 0} \\ \hline
0010010& 1 & 2 & {\color{red} 0} & {\color{red} 0} & {\color{green} 0} & {\color{green} 0} &{\color{blue} 0} &{\color{blue} 0} \\ \hline
0010011& 1 & 1 & {\color{red} 0} & {\color{red} 0} & {\color{green} 0} & {\color{green} 0} &{\color{blue} 1} &{\color{blue} 0} \\ \hline
0010100& 1 & 1 & {\color{red} 0} & {\color{red} 0} & {\color{green} 1} & {\color{green} 1} &{\color{blue} 0} &{\color{blue} 0} \\ \hline
0010101& 0 & 1 & {\color{red} 0} & {\color{red} 0} & {\color{green} 1} & {\color{green} 2} &{\color{blue} 0} &{\color{blue} 0} \\ \hline
0010110& 0 & 1 & {\color{red} 0} & {\color{red} 1} & {\color{green} 0} & {\color{green} 1} &{\color{blue} 0} &{\color{blue} 0} \\ \hline
0010111& 0 & 1 & {\color{red} 0} & {\color{red} 1} & {\color{green} 0} & {\color{green} 1} &{\color{blue} 0} &{\color{blue} 0} \\ \hline
0100000& 1 & 0 & {\color{red} 0} & {\color{red} 0} & {\color{green} 0} & {\color{green} 0} &{\color{blue} 0} &{\color{blue} 0} \\ \hline
0100001& 1 & 0 & {\color{red} 0} & {\color{red} 0} & {\color{green} 0} & {\color{green} 0} &{\color{blue} 0} &{\color{blue} 0} \\ \hline
0100010& 1 & 1 & {\color{red} 0} & {\color{red} 0} & {\color{green} 0} & {\color{green} 0} &{\color{blue} 0} &{\color{blue} 0} \\ \hline
0100011& 1 & 0 & {\color{red} 0} & {\color{red} 0} & {\color{green} 0} & {\color{green} 0} &{\color{blue} 1} &{\color{blue} 0} \\ \hline
0100100& 2 & 1 & {\color{red} 0} & {\color{red} 0} & {\color{green} 0} & {\color{green} 0} &{\color{blue} 0} &{\color{blue} 0} \\ \hline
0100101& 1 & 1 & {\color{red} 0} & {\color{red} 0} & {\color{green} 0} & {\color{green} 1} &{\color{blue} 0} &{\color{blue} 0} \\ \hline
0100110& 1 & 0 & {\color{red} 0} & {\color{red} 0} & {\color{green} 0} & {\color{green} 0} &{\color{blue} 1} &{\color{blue} 0} \\ \hline
0100111& 1 & 0 & {\color{red} 0} & {\color{red} 0} & {\color{green} 0} & {\color{green} 0} &{\color{blue} 1} &{\color{blue} 0} \\ \hline
0110000& 0 & 0 & {\color{red} 0} & {\color{red} 0} & {\color{green} 0} & {\color{green} 0} &{\color{blue} 0} &{\color{blue} 1} \\ \hline
0110001& 0 & 0 & {\color{red} 0} & {\color{red} 0} & {\color{green} 0} & {\color{green} 0} &{\color{blue} 0} &{\color{blue} 1} \\ \hline
0110010& 0 & 1 & {\color{red} 0} & {\color{red} 0} & {\color{green} 0} & {\color{green} 0} &{\color{blue} 0} &{\color{blue} 1} \\ \hline
0110011& 0 & 0 & {\color{red} 0} & {\color{red} 0} & {\color{green} 0} & {\color{green} 0} &{\color{blue} 1} &{\color{blue} 1} \\ \hline
0110100& 1 & 0 & {\color{red} 1} & {\color{red} 0} & {\color{green} 1} & {\color{green} 0} &{\color{blue} 0} &{\color{blue} 0} \\ \hline
0110101& 0 & 0 & {\color{red} 1} & {\color{red} 0} & {\color{green} 1} & {\color{green} 1} &{\color{blue} 0} &{\color{blue} 0} \\ \hline
0110110& 0 & 0 & {\color{red} 1} & {\color{red} 1} & {\color{green} 0} & {\color{green} 0} &{\color{blue} 0} &{\color{blue} 0} \\ \hline
0110111& 0 & 0 & {\color{red} 1} & {\color{red} 1} & {\color{green} 0} & {\color{green} 0} &{\color{blue} 0} &{\color{blue} 0} \\ \hline
1000000& 0 & 0 & {\color{red} 0} & {\color{red} 0} & {\color{green} 0} & {\color{green} 0} &{\color{blue} 0} &{\color{blue} 0} \\ \hline
1000001& 0 & 0 & {\color{red} 0} & {\color{red} 0} & {\color{green} 0} & {\color{green} 0} &{\color{blue} 0} &{\color{blue} 0} \\ \hline
1000010& 0 & 1 & {\color{red} 0} & {\color{red} 0} & {\color{green} 0} & {\color{green} 0} &{\color{blue} 0} &{\color{blue} 0} \\ \hline
1000011& 0 & 0 & {\color{red} 0} & {\color{red} 0} & {\color{green} 0} & {\color{green} 0} &{\color{blue} 1} &{\color{blue} 0} \\ \hline
1000100& 1 & 1 & {\color{red} 0} & {\color{red} 0} & {\color{green} 0} & {\color{green} 0} &{\color{blue} 0} &{\color{blue} 0} \\ \hline
1000101& 0 & 1 & {\color{red} 0} & {\color{red} 0} & {\color{green} 0} & {\color{green} 1} &{\color{blue} 0} &{\color{blue} 0} \\ \hline
1000110& 0 & 0 & {\color{red} 0} & {\color{red} 0} & {\color{green} 0} & {\color{green} 0} &{\color{blue} 1} &{\color{blue} 0} \\ \hline
1000111& 0 & 0 & {\color{red} 0} & {\color{red} 0} & {\color{green} 0} & {\color{green} 0} &{\color{blue} 1} &{\color{blue} 0} \\ \hline
1010000& 1 & 0 & {\color{red} 0} & {\color{red} 0} & {\color{green} 1} & {\color{green} 0} &{\color{blue} 0} &{\color{blue} 0} \\ \hline
1010001& 1 & 0 & {\color{red} 0} & {\color{red} 0} & {\color{green} 1} & {\color{green} 0} &{\color{blue} 0} &{\color{blue} 0} \\ \hline
1010010& 1 & 1 & {\color{red} 0} & {\color{red} 0} & {\color{green} 1} & {\color{green} 0} &{\color{blue} 0} &{\color{blue} 0} \\ \hline
1010011& 1 & 0 & {\color{red} 0} & {\color{red} 0} & {\color{green} 1} & {\color{green} 0} &{\color{blue} 1} &{\color{blue} 0} \\ \hline
1010100& 1 & 0 & {\color{red} 0} & {\color{red} 0} & {\color{green} 2} & {\color{green} 1} &{\color{blue} 0} &{\color{blue} 0} \\ \hline
1010101& 0 & 0 & {\color{red} 0} & {\color{red} 0} & {\color{green} 2} & {\color{green} 2} &{\color{blue} 0} &{\color{blue} 0} \\ \hline
1010110& 0 & 0 & {\color{red} 0} & {\color{red} 1} & {\color{green} 1} & {\color{green} 1} &{\color{blue} 0} &{\color{blue} 0} \\ \hline
1010111& 0 & 0 & {\color{red} 0} & {\color{red} 1} & {\color{green} 1} & {\color{green} 1} &{\color{blue} 0} &{\color{blue} 0} \\ \hline
1100000& 0 & 0 & {\color{red} 0} & {\color{red} 0} & {\color{green} 0} & {\color{green} 0} &{\color{blue} 0} &{\color{blue} 1} \\ \hline
1100001& 0 & 0 & {\color{red} 0} & {\color{red} 0} & {\color{green} 0} & {\color{green} 0} &{\color{blue} 0} &{\color{blue} 1} \\ \hline
1100010& 0 & 1 & {\color{red} 0} & {\color{red} 0} & {\color{green} 0} & {\color{green} 0} &{\color{blue} 0} &{\color{blue} 1} \\ \hline
1100011& 0 & 0 & {\color{red} 0} & {\color{red} 0} & {\color{green} 0} & {\color{green} 0} &{\color{blue} 1} &{\color{blue} 1} \\ \hline
1100100& 1 & 1 & {\color{red} 0} & {\color{red} 0} & {\color{green} 0} & {\color{green} 0} &{\color{blue} 0} &{\color{blue} 1} \\ \hline
1100101& 0 & 1 & {\color{red} 0} & {\color{red} 0} & {\color{green} 0} & {\color{green} 1} &{\color{blue} 0} &{\color{blue} 1} \\ \hline
1100110& 0 & 0 & {\color{red} 0} & {\color{red} 0} & {\color{green} 0} & {\color{green} 0} &{\color{blue} 1} &{\color{blue} 1} \\ \hline
1100111& 0 & 0 & {\color{red} 0} & {\color{red} 0} & {\color{green} 0} & {\color{green} 0} &{\color{blue} 0} &{\color{blue} 0} \\ \hline
1110000& 0 & 0 & {\color{red} 0} & {\color{red} 0} & {\color{green} 0} & {\color{green} 0} &{\color{blue} 0} &{\color{blue} 1} \\ \hline
1110001& 0 & 0 & {\color{red} 0} & {\color{red} 0} & {\color{green} 0} & {\color{green} 0} &{\color{blue} 0} &{\color{blue} 1} \\ \hline
1110010& 0 & 1 & {\color{red} 0} & {\color{red} 0} & {\color{green} 0} & {\color{green} 0} &{\color{blue} 0} &{\color{blue} 1} \\ \hline
1110011& 0 & 0 & {\color{red} 0} & {\color{red} 0} & {\color{green} 0} & {\color{green} 0} &{\color{blue} 1} &{\color{blue} 1} \\ \hline
1110100& 1 & 0 & {\color{red} 1} & {\color{red} 0} & {\color{green} 1} & {\color{green} 0} &{\color{blue} 0} &{\color{blue} 0} \\ \hline
1110101& 0 & 0 & {\color{red} 1} & {\color{red} 0} & {\color{green} 1} & {\color{green} 1} &{\color{blue} 0} &{\color{blue} 0} \\ \hline
1110110& 0 & 0 & {\color{red} 1} & {\color{red} 1} & {\color{green} 0} & {\color{green} 0} &{\color{blue} 0} &{\color{blue} 0} \\ \hline
1110111& 0 & 0 & {\color{red} 1} & {\color{red} 1} & {\color{green} 0} & {\color{green} 0} &{\color{blue} 0} &{\color{blue} 0} \\ \hline

	\caption{Sequences adjoined, $d_n=0$}
	\label{tab:sadj0}
\end{longtable}

The table \ref{tab:sadj1} shows the number of  sequences adjoined for each case, when $d_n=1$ ($P'' = d_{n-3},\,d_{n-2},\,d_{n-1},\,1,\,d_0,\,d_1,\,d_2$).

	\begin{longtable}{|c|c|c|c|c|c|c|c|c|}			
\hline
	$P''$	& 0100 & 0010 & {\color{red}1101} & {\color{red}1011}& {\color{green} 1010} & {\color{green}0101} & {\color{blue} 0011} & {\color{blue} 1100} \\
		\hline
		0001000 & 1 & 1 & {\color{red} 0} & {\color{red} 0} & {\color{green} 0} & {\color{green} 0} &{\color{blue} 0} &{\color{blue} 0} \\ \hline
	0001001	& 1 & 1 & {\color{red} 0} & {\color{red} 0} & {\color{green} 0} & {\color{green} 0} &{\color{blue} 0} &{\color{blue} 0} \\ \hline
0001010& 0 & 1 & {\color{red} 0} & {\color{red} 0} & {\color{green} 1} & {\color{green} 1} &{\color{blue} 0} &{\color{blue} 0} \\ \hline
0001011& 0 & 1 & {\color{red} 0} & {\color{red} 1} & {\color{green} 0} & {\color{green} 1} &{\color{blue} 0} &{\color{blue} 0} \\ \hline
0001100& 0 & 0 & {\color{red} 0} & {\color{red} 0} & {\color{green} 0} & {\color{green} 0} &{\color{blue} 1} &{\color{blue} 1} \\ \hline
0001101& 0 & 0 & {\color{red} 1} & {\color{red} 0} & {\color{green} 0} & {\color{green} 0} &{\color{blue} 1} &{\color{blue} 0} \\ \hline
0001110& 0 & 0 & {\color{red} 0} & {\color{red} 0} & {\color{green} 0} & {\color{green} 0} &{\color{blue} 1} &{\color{blue} 0} \\ \hline
0001111& 0 & 0 & {\color{red} 0} & {\color{red} 0} & {\color{green} 0} & {\color{green} 0} &{\color{blue} 1} &{\color{blue} 0} \\ \hline
0011000& 0 & 0 & {\color{red} 0} & {\color{red} 0} & {\color{green} 0} & {\color{green} 0} &{\color{blue} 1} &{\color{blue} 1} \\ \hline
0011001& 0 & 0 & {\color{red} 0} & {\color{red} 0} & {\color{green} 0} & {\color{green} 0} &{\color{blue} 1} &{\color{blue} 1} \\ \hline
0011010& 0 & 0 & {\color{red} 1} & {\color{red} 0} & {\color{green} 1} & {\color{green} 0} &{\color{blue} 1} &{\color{blue} 0} \\ \hline
0011011& 0 & 0 & {\color{red} 1} & {\color{red} 1} & {\color{green} 0} & {\color{green} 0} &{\color{blue} 1} &{\color{blue} 0} \\ \hline
0011100& 0 & 0 & {\color{red} 0} & {\color{red} 0} & {\color{green} 0} & {\color{green} 0} &{\color{blue} 1} &{\color{blue} 1} \\ \hline
0011101& 0 & 0 & {\color{red} 1} & {\color{red} 0} & {\color{green} 0} & {\color{green} 0} &{\color{blue} 1} &{\color{blue} 0} \\ \hline
0011110& 0 & 0 & {\color{red} 0} & {\color{red} 0} & {\color{green} 0} & {\color{green} 0} &{\color{blue} 1} &{\color{blue} 0} \\ \hline
0011111& 0 & 0 & {\color{red} 0} & {\color{red} 0} & {\color{green} 0} & {\color{green} 0} &{\color{blue} 1} &{\color{blue} 0} \\ \hline
0101000& 1 & 0 & {\color{red} 0} & {\color{red} 0} & {\color{green} 1} & {\color{green} 1} &{\color{blue} 0} &{\color{blue} 0} \\ \hline
0101001& 1 & 0 & {\color{red} 0} & {\color{red} 0} & {\color{green} 1} & {\color{green} 1} &{\color{blue} 0} &{\color{blue} 0} \\ \hline
0101010& 0 & 0 & {\color{red} 0} & {\color{red} 0} & {\color{green} 2} & {\color{green} 2} &{\color{blue} 0} &{\color{blue} 0} \\ \hline
0101011& 0 & 0 & {\color{red} 0} & {\color{red} 1} & {\color{green} 1} & {\color{green} 2} &{\color{blue} 0} &{\color{blue} 0} \\ \hline
0101100& 0 & 0 & {\color{red} 0} & {\color{red} 1} & {\color{green} 0} & {\color{green} 1} &{\color{blue} 0} &{\color{blue} 1} \\ \hline
0101101& 0 & 0 & {\color{red} 1} & {\color{red} 1} & {\color{green} 0} & {\color{green} 1} &{\color{blue} 0} &{\color{blue} 0} \\ \hline
0101110& 0 & 0 & {\color{red} 0} & {\color{red} 1} & {\color{green} 0} & {\color{green} 1} &{\color{blue} 0} &{\color{blue} 0} \\ \hline
0101111& 0 & 0 & {\color{red} 0} & {\color{red} 1} & {\color{green} 0} & {\color{green} 1} &{\color{blue} 0} &{\color{blue} 0} \\ \hline
0111000& 0 & 0 & {\color{red} 0} & {\color{red} 0} & {\color{green} 0} & {\color{green} 0} &{\color{blue} 0} &{\color{blue} 1} \\ \hline
0111001& 0 & 0 & {\color{red} 0} & {\color{red} 0} & {\color{green} 0} & {\color{green} 0} &{\color{blue} 0} &{\color{blue} 1} \\ \hline
0111010& 0 & 0 & {\color{red} 1} & {\color{red} 0} & {\color{green} 1} & {\color{green} 0} &{\color{blue} 0} &{\color{blue} 0} \\ \hline
0111011& 0 & 0 & {\color{red} 1} & {\color{red} 1} & {\color{green} 0} & {\color{green} 0} &{\color{blue} 0} &{\color{blue} 0} \\ \hline
0111100& 0 & 0 & {\color{red} 0} & {\color{red} 0} & {\color{green} 0} & {\color{green} 0} &{\color{blue} 0} &{\color{blue} 1} \\ \hline
0111101& 0 & 0 & {\color{red} 1} & {\color{red} 0} & {\color{green} 0} & {\color{green} 0} &{\color{blue} 0} &{\color{blue} 0} \\ \hline
0111110& 0 & 0 & {\color{red} 0} & {\color{red} 0} & {\color{green} 0} & {\color{green} 0} &{\color{blue} 0} &{\color{blue} 0} \\ \hline
0111111& 0 & 0 & {\color{red} 0} & {\color{red} 0} & {\color{green} 0} & {\color{green} 0} &{\color{blue} 0} &{\color{blue} 0} \\ \hline
1001000& 1 & 1 & {\color{red} 0} & {\color{red} 0} & {\color{green} 0} & {\color{green} 0} &{\color{blue} 0} &{\color{blue} 0} \\ \hline
1001001& 1 & 1 & {\color{red} 0} & {\color{red} 0} & {\color{green} 0} & {\color{green} 0} &{\color{blue} 0} &{\color{blue} 0} \\ \hline
1001010& 0 & 1 & {\color{red} 0} & {\color{red} 0} & {\color{green} 1} & {\color{green} 1} &{\color{blue} 0} &{\color{blue} 0} \\ \hline
1001011& 0 & 1 & {\color{red} 0} & {\color{red} 1} & {\color{green} 0} & {\color{green} 1} &{\color{blue} 0} &{\color{blue} 0} \\ \hline
1001100& 1 & 1 & {\color{red} 0} & {\color{red} 0} & {\color{green} 0} & {\color{green} 0} &{\color{blue} 1} &{\color{blue} 1} \\ \hline
1001101& 0 & 0 & {\color{red} 1} & {\color{red} 0} & {\color{green} 0} & {\color{green} 0} &{\color{blue} 1} &{\color{blue} 0} \\ \hline
1001110& 0 & 0 & {\color{red} 0} & {\color{red} 0} & {\color{green} 0} & {\color{green} 0} &{\color{blue} 1} &{\color{blue} 0} \\ \hline
1001111& 0 & 0 & {\color{red} 0} & {\color{red} 0} & {\color{green} 0} & {\color{green} 0} &{\color{blue} 1} &{\color{blue} 0} \\ \hline
1011000& 0 & 0 & {\color{red} 0} & {\color{red} 1} & {\color{green} 0} & {\color{green} 0} &{\color{blue} 0} &{\color{blue} 1} \\ \hline
1011001& 0 & 0 & {\color{red} 0} & {\color{red} 1} & {\color{green} 0} & {\color{green} 0} &{\color{blue} 0} &{\color{blue} 1} \\ \hline
1011010& 0 & 0 & {\color{red} 1} & {\color{red} 1} & {\color{green} 1} & {\color{green} 0} &{\color{blue} 0} &{\color{blue} 0} \\ \hline
1011011& 0 & 0 & {\color{red} 1} & {\color{red} 2} & {\color{green} 0} & {\color{green} 0} &{\color{blue} 0} &{\color{blue} 0} \\ \hline
1011100& 0 & 0 & {\color{red} 0} & {\color{red} 1} & {\color{green} 0} & {\color{green} 0} &{\color{blue} 0} &{\color{blue} 1} \\ \hline
1011101& 0 & 0 & {\color{red} 1} & {\color{red} 1} & {\color{green} 0} & {\color{green} 0} &{\color{blue} 0} &{\color{blue} 0} \\ \hline
1011110& 0 & 0 & {\color{red} 0} & {\color{red} 1} & {\color{green} 0} & {\color{green} 0} &{\color{blue} 0} &{\color{blue} 0} \\ \hline
1011111& 0 & 0 & {\color{red} 0} & {\color{red} 1} & {\color{green} 0} & {\color{green} 0} &{\color{blue} 0} &{\color{blue} 0} \\ \hline
1101000& 1 & 0 & {\color{red} 1} & {\color{red} 0} & {\color{green} 1} & {\color{green} 0} &{\color{blue} 0} &{\color{blue} 0} \\ \hline
1101001& 1 & 0 & {\color{red} 1} & {\color{red} 0} & {\color{green} 1} & {\color{green} 0} &{\color{blue} 0} &{\color{blue} 0} \\ \hline
1101010& 0 & 0 & {\color{red} 1} & {\color{red} 0} & {\color{green} 2} & {\color{green} 1} &{\color{blue} 0} &{\color{blue} 0} \\ \hline
1101011& 0 & 0 & {\color{red} 1} & {\color{red} 1} & {\color{green} 1} & {\color{green} 1} &{\color{blue} 0} &{\color{blue} 0} \\ \hline
1101100& 0 & 0 & {\color{red} 1} & {\color{red} 1} & {\color{green} 0} & {\color{green} 0} &{\color{blue} 0} &{\color{blue} 1} \\ \hline
1101101& 0 & 0 & {\color{red} 2} & {\color{red} 1} & {\color{green} 0} & {\color{green} 0} &{\color{blue} 0} &{\color{blue} 0} \\ \hline
1101110& 0 & 0 & {\color{red} 1} & {\color{red} 1} & {\color{green} 0} & {\color{green} 0} &{\color{blue} 0} &{\color{blue} 0} \\ \hline
1101111& 0 & 0 & {\color{red} 1} & {\color{red} 1} & {\color{green} 0} & {\color{green} 0} &{\color{blue} 0} &{\color{blue} 0} \\ \hline
1111000& 0 & 0 & {\color{red} 0} & {\color{red} 0} & {\color{green} 0} & {\color{green} 0} &{\color{blue} 0} &{\color{blue} 1} \\ \hline
1111001& 0 & 0 & {\color{red} 0} & {\color{red} 0} & {\color{green} 0} & {\color{green} 0} &{\color{blue} 0} &{\color{blue} 1} \\ \hline
1111010& 0 & 0 & {\color{red} 1} & {\color{red} 0} & {\color{green} 1} & {\color{green} 0} &{\color{blue} 0} &{\color{blue} 0} \\ \hline
1111011& 0 & 0 & {\color{red} 1} & {\color{red} 1} & {\color{green} 0} & {\color{green} 0} &{\color{blue} 0} &{\color{blue} 0} \\ \hline
1111100& 0 & 0 & {\color{red} 0} & {\color{red} 0} & {\color{green} 0} & {\color{green} 0} &{\color{blue} 0} &{\color{blue} 1} \\ \hline
1111101& 0 & 0 & {\color{red} 1} & {\color{red} 0} & {\color{green} 0} & {\color{green} 0} &{\color{blue} 0} &{\color{blue} 0} \\ \hline
1111110& 0 & 0 & {\color{red} 0} & {\color{red} 0} & {\color{green} 0} & {\color{green} 0} &{\color{blue} 0} &{\color{blue} 0} \\ \hline
1111111& 0 & 0 & {\color{red} 0} & {\color{red} 0} & {\color{green} 0} & {\color{green} 0} &{\color{blue} 0} &{\color{blue} 0} \\ \hline

	\caption{Sequences adjoined, $d_n=1$}
	\label{tab:sadj1}
\end{longtable}

Now, there is only to prove that the difference between the number
of occurences of the two terms for each periodic sequence is the same
for every  $P'$  and $P''$. For each case, there are the
sequences lost to remove (computed in table \ref{tab:spd}) and
the sequences adjoined (computed in table \ref{tab:sadj0} when $d_n =0$
and in table \ref{tab:sadj1} when $d_n=1$). The results
are displayed in table \ref{tab:ecaj0} and \ref{tab:ecaj1}. In order to understand better these tables, consider as example $P'=1010101$ in table
\ref{tab:ecaj0}. By moving from 
$D=101,\,d_3,\,d_4,\,\dots\,d_{n-4},\,101$ to                       $D'=101,\,d_3,\,d_4,\,\dots\,d_{n-4},\,1010$  we lose one sequence
$1101$ and one $1011$, as shown by table \ref{tab:spd}
($P=101101$). Nevertheless, we gain two sequences  $1010$  and two sequences
$0101$, as shown by table \ref{tab:ecaj0}
($P'=1010101$). In this case, the difference between the sequences is the same for $D$ and $D'$; if there was a $\Delta$Difference of  +1, this means that if the sequence $D$ has a difference $E$ between the number of occurences of
the two terms for each periodic sequence, the sequence $D'$ has a difference
of $E+1$.

	\begin{longtable}{|c|c|c|c|c|c|c|c|c|c|}			
\hline
	$P'$	& 0100 & 0010 & {\color{red}1101} & {\color{red}1011}& {\color{green} 1010} & {\color{green}0101} & {\color{blue} 0011} & {\color{blue} 1100}& $\Delta$Difference \\
		\hline
		0000000 & 0 & 0 & {\color{red} 0} & {\color{red} 0} & {\color{green} 0} & {\color{green} 0} &{\color{blue} 0} &{\color{blue} 0} & 0\\ \hline
	0000001	& 0 & 0 & {\color{red} 0} & {\color{red} 0} & {\color{green} 0} & {\color{green} 0} &{\color{blue} 0} &{\color{blue} 0} & 0 \\ \hline
0000010& 0 & 1-1 & {\color{red} 0} & {\color{red} 0} & {\color{green} 0} & {\color{green} 0} &{\color{blue} 0} &{\color{blue} 0} & 0 \\ \hline
0000011& 0 & 0 & {\color{red} 0} & {\color{red} 0} & {\color{green} 0} & {\color{green} 0} &{\color{blue} 1-1} &{\color{blue} 0} &0 \\ \hline
0000100& 1-1 & 1-1 & {\color{red} 0} & {\color{red} 0} & {\color{green} 0} & {\color{green} 0} &{\color{blue} 0} &{\color{blue} 0} &0 \\ \hline
0000101& 0 & 1-1 & {\color{red} 0} & {\color{red} 0} & {\color{green} 0} & {\color{green} 1-1} &{\color{blue} 0} &{\color{blue} 0} &0 \\ \hline
0000110& 0 & 0 & {\color{red} 0} & {\color{red} 0} & {\color{green} 0} & {\color{green} 0} &{\color{blue} 1-1} &{\color{blue} 0} &0\\ \hline
0000111& 0 & 0 & {\color{red} 0} & {\color{red} 0} & {\color{green} 0} & {\color{green} 0} &{\color{blue} 1-1} &{\color{blue} 0} &0\\ \hline
0010000& 1-1 & 1-1 & {\color{red} 0} & {\color{red} 0} & {\color{green} 0} & {\color{green} 0} &{\color{blue} 0} &{\color{blue} 0} &0\\ \hline
0010001& 1-1 & 1-1 & {\color{red} 0} & {\color{red} 0} & {\color{green} 0} & {\color{green} 0} &{\color{blue} 0} &{\color{blue} 0} &0\\ \hline
0010010& 1 & 2-1 & {\color{red} 0} & {\color{red} 0} & {\color{green} -1} & {\color{green} -1} &{\color{blue} 0} &{\color{blue} 0}  & 0\\ \hline
0010011& 1 & 1-1 & {\color{red} 0} & {\color{red} -1} & {\color{green} 0} & {\color{green} -1} &{\color{blue} 1} &{\color{blue} 0} & +1 \\ \hline
0010100& 1 & 1 & {\color{red} 0} & {\color{red} 0} & {\color{green} 1} & {\color{green} 1} &{\color{blue} -1} &{\color{blue} -1} & 0\\ \hline
0010101& 0 & 1 & {\color{red} -1} & {\color{red} 0} & {\color{green} 1} & {\color{green} 2} &{\color{blue} -1} &{\color{blue} 0} & -1 \\ \hline
0010110& 0 & 1 & {\color{red} 0} & {\color{red} 1} & {\color{green} 0} & {\color{green} 1} &{\color{blue} -1} &{\color{blue} 0} & -1\\ \hline
0010111& 0 & 1 & {\color{red} 0} & {\color{red} 1} & {\color{green} 0} & {\color{green} 1} &{\color{blue} -1} &{\color{blue} 0} & -1 \\ \hline
0100000& 1-1 & 0 & {\color{red} 0} & {\color{red} 0} & {\color{green} 0} & {\color{green} 0} &{\color{blue} 0} &{\color{blue} 0} & 0 \\ \hline
0100001& 1-1 & 0 & {\color{red} 0} & {\color{red} 0} & {\color{green} 0} & {\color{green} 0} &{\color{blue} 0} &{\color{blue} 0} & 0 \\ \hline
0100010& 1-1 & 1-1 & {\color{red} 0} & {\color{red} 0} & {\color{green} 0} & {\color{green} 0} &{\color{blue} 0} &{\color{blue} 0} & 0 \\ \hline
0100011& 1-1 & 0 & {\color{red} 0} & {\color{red} 0} & {\color{green} 0} & {\color{green} 0} &{\color{blue} 1-1} &{\color{blue} 0} & 0 \\ \hline
0100100& 2-1 & 1 & {\color{red} 0} & {\color{red} 0} & {\color{green} -1} & {\color{green} -1} &{\color{blue} 0} &{\color{blue} 0} & 0\\ \hline
0100101& 1 & 1 & {\color{red} 0} & {\color{red} 0} & {\color{green} -1} & {\color{green} 1-2} &{\color{blue} 0} &{\color{blue} 0} & 0\\ \hline
0100110& 1 & 0 & {\color{red} 0} & {\color{red} -1} & {\color{green} 0} & {\color{green} -1} &{\color{blue} 1} &{\color{blue} 0} & +1\\ \hline
0100111& 1 & 0 & {\color{red} 0} & {\color{red} -1} & {\color{green} 0} & {\color{green} -1} &{\color{blue} 1} &{\color{blue} 0} & +1\\ \hline
0110000& 0 & 0 & {\color{red} 0} & {\color{red} 0} & {\color{green} 0} & {\color{green} 0} &{\color{blue} 0} &{\color{blue} 1-1} & 0 \\ \hline
0110001& 0 & 0 & {\color{red} 0} & {\color{red} 0} & {\color{green} 0} & {\color{green} 0} &{\color{blue} 0} &{\color{blue} 1-1} & 0 \\ \hline
0110010& 0 & 1 & {\color{red} -1} & {\color{red} 0} & {\color{green} -1} & {\color{green} 0} &{\color{blue} 0} &{\color{blue} 1} & -1 \\ \hline
0110011& 0 & 0 & {\color{red} -1} & {\color{red} -1} & {\color{green} 0} & {\color{green} 0} &{\color{blue} 1} &{\color{blue} 1} & 0\\ \hline
0110100& 1 & 0 & {\color{red} 1} & {\color{red} 0} & {\color{green} 1} & {\color{green} 0} &{\color{blue} 0} &{\color{blue} -1} & 1\\ \hline
0110101& 0 & 0 & {\color{red} 1-1} & {\color{red} 0} & {\color{green} 1} & {\color{green} 1} &{\color{blue} 0} &{\color{blue} 0} & 0\\ \hline
0110110& 0 & 0 & {\color{red} 1} & {\color{red} 1} & {\color{green} 0} & {\color{green} 0} &{\color{blue} 0} &{\color{blue} 0} & 0\\ \hline
0110111& 0 & 0 & {\color{red} 1} & {\color{red} 1} & {\color{green} 0} & {\color{green} 0} &{\color{blue} 0} &{\color{blue} 0} & 0\\ \hline
1000000& 0 & 0 & {\color{red} 0} & {\color{red} 0} & {\color{green} 0} & {\color{green} 0} &{\color{blue} 0} &{\color{blue} 0} & 0\\ \hline
1000001& 0 & 0 & {\color{red} 0} & {\color{red} 0} & {\color{green} 0} & {\color{green} 0} &{\color{blue} 0} &{\color{blue} 0} & 0\\ \hline
1000010& 0 & 1-1 & {\color{red} 0} & {\color{red} 0} & {\color{green} 0} & {\color{green} 0} &{\color{blue} 0} &{\color{blue} 0} & 0 \\ \hline
1000011& 0 & 0 & {\color{red} 0} & {\color{red} 0} & {\color{green} 0} & {\color{green} 0} &{\color{blue} 1-1} &{\color{blue} 0} & 0 \\ \hline
1000100& 1-1 & 1-1 & {\color{red} 0} & {\color{red} 0} & {\color{green} 0} & {\color{green} 0} &{\color{blue} 0} &{\color{blue} 0} & 0 \\ \hline
1000101& 0 & 1-1 & {\color{red} 0} & {\color{red} 0} & {\color{green} 0} & {\color{green} 1-1} &{\color{blue} 0} &{\color{blue} 0} &0\\ \hline
1000110& 0 & 0 & {\color{red} 0} & {\color{red} 0} & {\color{green} 0} & {\color{green} 0} &{\color{blue} 1-1} &{\color{blue} 0} & 0\\ \hline
1000111& 0 & 0 & {\color{red} 0} & {\color{red} 0} & {\color{green} 0} & {\color{green} 0} &{\color{blue} 1-1} &{\color{blue} 0} & 0\\ \hline
1010000& 1-1 & 0 & {\color{red} 0} & {\color{red} 0} & {\color{green} 1-1} & {\color{green} 0} &{\color{blue} 0} &{\color{blue} 0} & 0\\ \hline
1010001& 1-1 & 0 & {\color{red} 0} & {\color{red} 0} & {\color{green} 1-1} & {\color{green} 0} &{\color{blue} 0} &{\color{blue} 0} & 0\\ \hline
1010010& 1 & 1 & {\color{red} 0} & {\color{red} 0} & {\color{green} 1-2} & {\color{green} -1} &{\color{blue} 0} &{\color{blue} 0} & 0\\ \hline
1010011& 1 & 0 & {\color{red} 0} & {\color{red} -1} & {\color{green} 1-1} & {\color{green} -1} &{\color{blue} 1} &{\color{blue} 0} & +1\\ \hline
1010100& 1 & 0 & {\color{red} 0} & {\color{red} -1} & {\color{green} 2} & {\color{green} 1} &{\color{blue} 0} &{\color{blue} -1} & +1\\ \hline
1010101& 0 & 0 & {\color{red} -1} & {\color{red} -1} & {\color{green} 2} & {\color{green} 2} &{\color{blue} 0} &{\color{blue} 0} & 0\\ \hline
1010110& 0 & 0 & {\color{red} 0} & {\color{red} 1-1} & {\color{green} 1} & {\color{green} 1} &{\color{blue} 0} &{\color{blue} 0}& 0 \\ \hline
1010111& 0 & 0 & {\color{red} 0} & {\color{red} 1-1} & {\color{green} 1} & {\color{green} 1} &{\color{blue} 0} &{\color{blue} 0} & 0\\ \hline
1100000& 0 & 0 & {\color{red} 0} & {\color{red} 0} & {\color{green} 0} & {\color{green} 0} &{\color{blue} 0} &{\color{blue} 1-1} & 0\\ \hline
1100001& 0 & 0 & {\color{red} 0} & {\color{red} 0} & {\color{green} 0} & {\color{green} 0} &{\color{blue} 0} &{\color{blue} 1-1}  &0\\ \hline
1100010& 0 & 1-1 & {\color{red} 0} & {\color{red} 0} & {\color{green} 0} & {\color{green} 0} &{\color{blue} 0} &{\color{blue} 1-1} & 0 \\ \hline
1100011& 0 & 0 & {\color{red} 0} & {\color{red} 0} & {\color{green} 0} & {\color{green} 0} &{\color{blue} 1-1} &{\color{blue} 1-1}  &0\\ \hline
1100100& 1-1 & 1 & {\color{red} -1} & {\color{red} 0} & {\color{green} -1} & {\color{green} 0} &{\color{blue} 0} &{\color{blue} 1} &-1\\ \hline
1100101& 0 & 1 & {\color{red} -1} & {\color{red} 0} & {\color{green} -1} & {\color{green} 1-1} &{\color{blue} 0} &{\color{blue} 1} & -1\\ \hline
1100110& 0 & 0 & {\color{red} -1} & {\color{red} -1} & {\color{green} 0} & {\color{green} 0} &{\color{blue} 1} &{\color{blue} 1} & 0 \\ \hline
1100111& 0 & 0 & {\color{red} -1} & {\color{red} -1} & {\color{green} 0} & {\color{green} 0} &{\color{blue} 0} &{\color{blue} 0} & 0\\ \hline
1110000& 0 & 0 & {\color{red} 0} & {\color{red} 0} & {\color{green} 0} & {\color{green} 0} &{\color{blue} 0} &{\color{blue} 1-1} & 0\\ \hline
1110001& 0 & 0 & {\color{red} 0} & {\color{red} 0} & {\color{green} 0} & {\color{green} 0} &{\color{blue} 0} &{\color{blue} 1-1} & 0 \\ \hline
1110010& 0 & 1 & {\color{red} -1} & {\color{red} 0} & {\color{green} -1} & {\color{green} 0} &{\color{blue} 0} &{\color{blue} 1} & 0\\ \hline
1110011& 0 & 0 & {\color{red} -1} & {\color{red} -1} & {\color{green} 0} & {\color{green} 0} &{\color{blue} 1} &{\color{blue} 1} & 0\\ \hline
1110100& 1 & 0 & {\color{red} 1} & {\color{red} 0} & {\color{green} 1} & {\color{green} 0} &{\color{blue} 0} &{\color{blue} -1} & +1 \\ \hline
1110101& 0 & 0 & {\color{red} 1-1} & {\color{red} 0} & {\color{green} 1} & {\color{green} 1} &{\color{blue} 0} &{\color{blue} 0} & 0\\ \hline
1110110& 0 & 0 & {\color{red} 1} & {\color{red} 1} & {\color{green} 0} & {\color{green} 0} &{\color{blue} 0} &{\color{blue} 0} & 0\\ \hline
1110111& 0 & 0 & {\color{red} 1} & {\color{red} 1} & {\color{green} 0} & {\color{green} 0} &{\color{blue} 0} &{\color{blue} 0} & 0\\ \hline

	\caption{$\Delta$Difference, $d_n=0$}
	\label{tab:ecaj0}
\end{longtable}

	\begin{longtable}{|c|c|c|c|c|c|c|c|c|c|}			
\hline
	$P''$	& 0100 & 0010 & {\color{red}1101} & {\color{red}1011}& {\color{green} 1010} & {\color{green}0101} & {\color{blue} 0011} & {\color{blue} 1100} & $\Delta$Difference \\
		\hline
		0001000 & 1 & 1 & {\color{red} 0} & {\color{red} 0} & {\color{green} 0} & {\color{green} 0} &{\color{blue} 0} &{\color{blue} 0} & 0 \\ \hline
	0001001	& 1 & 1 & {\color{red} 0} & {\color{red} 0} & {\color{green} 0} & {\color{green} 0} &{\color{blue} 0} &{\color{blue} 0} & 0\\ \hline
0001010& 0 & 1-1 & {\color{red} 0} & {\color{red} 0} & {\color{green} 1} & {\color{green} 1} &{\color{blue} 0} &{\color{blue} 0} & 0\\ \hline
0001011& 0 & 1 & {\color{red} 0} & {\color{red} 1} & {\color{green} 0} & {\color{green} 1} &{\color{blue} -1} &{\color{blue} 0} & -1 \\ \hline
0001100& -1 & -1 & {\color{red} 0} & {\color{red} 0} & {\color{green} 0} & {\color{green} 0} &{\color{blue} 1} &{\color{blue} 1} & 0 \\ \hline
0001101& 0 & -1 & {\color{red} 1} & {\color{red} 0} & {\color{green} 0} & {\color{green} -1} &{\color{blue} 1} &{\color{blue} 0} & 1 \\ \hline
0001110& 0 & 0 & {\color{red} 0} & {\color{red} 0} & {\color{green} 0} & {\color{green} 0} &{\color{blue} 1-1} &{\color{blue} 0} & 0\\ \hline
0001111& 0 & 0 & {\color{red} 0} & {\color{red} 0} & {\color{green} 0} & {\color{green} 0} &{\color{blue} 1-1} &{\color{blue} 0} & 0 \\ \hline
0011000& -1 & -1 & {\color{red} 0} & {\color{red} 0} & {\color{green} 0} & {\color{green} 0} &{\color{blue} 1} &{\color{blue} 1} & 0 \\ \hline
0011001& -1 & -1 & {\color{red} 0} & {\color{red} 0} & {\color{green} 0} & {\color{green} 0} &{\color{blue} 1} &{\color{blue} 1} & 0 \\ \hline
0011010& 0 & -1 & {\color{red} 1} & {\color{red} 0} & {\color{green} 1-1} & {\color{green} -1} &{\color{blue} 1} &{\color{blue} 0} & +1 \\ \hline
0011011& 0 & -1 & {\color{red} 1} & {\color{red} 1-1} & {\color{green} 0} & {\color{green} -1} &{\color{blue} 1} &{\color{blue} 0} & +1\\ \hline
0011100& 0 & 0 & {\color{red} 0} & {\color{red} 0} & {\color{green} 0} & {\color{green} 0} &{\color{blue} 1-1} &{\color{blue} 1-1} & 0\\ \hline
0011101& 0 & 0 & {\color{red} 1-1} & {\color{red} 0} & {\color{green} 0} & {\color{green} 0} &{\color{blue} 1-1} &{\color{blue} 0} & 0 \\ \hline
0011110& 0 & 0 & {\color{red} 0} & {\color{red} 0} & {\color{green} 0} & {\color{green} 0} &{\color{blue} 1-1} &{\color{blue} 0} & 0 \\ \hline
0011111& 0 & 0 & {\color{red} 0} & {\color{red} 0} & {\color{green} 0} & {\color{green} 0} &{\color{blue} 1-1} &{\color{blue} 0}& 0 \\ \hline
0101000& 1-1 & 0 & {\color{red} 0} & {\color{red} 0} & {\color{green} 1} & {\color{green} 1} &{\color{blue} 0} &{\color{blue} 0} & 0 \\ \hline
0101001& 1-1 & 0 & {\color{red} 0} & {\color{red} 0} & {\color{green} 1} & {\color{green} 1} &{\color{blue} 0} &{\color{blue} 0} &0 \\ \hline
0101010& -1 & -1 & {\color{red} 0} & {\color{red} 0} & {\color{green} 2} & {\color{green} 2} &{\color{blue} 0} &{\color{blue} 0} & 0 \\ \hline
0101011& -1 & 0 & {\color{red} 0} & {\color{red} 1} & {\color{green} 1} & {\color{green} 2} &{\color{blue} -1} &{\color{blue} 0} & -1\\ \hline
0101100& -1 & 0 & {\color{red} 0} & {\color{red} 1} & {\color{green} -1} & {\color{green} 1-1} &{\color{blue} 0} &{\color{blue} 1} & -1\\ \hline
0101101& 0 & 0 & {\color{red} 1} & {\color{red} 1} & {\color{green} -1} & {\color{green} 1-2} &{\color{blue} 0} &{\color{blue} 0} & 0\\ \hline
0101110& 0 & 0 & {\color{red} 0} & {\color{red} 1-1} & {\color{green} 0} & {\color{green} 1-1} &{\color{blue} 0} &{\color{blue} 0} & 0\\ \hline
0101111& 0 & 0 & {\color{red} 0} & {\color{red} 1-1} & {\color{green} 0} & {\color{green} 1-1} &{\color{blue} 0} &{\color{blue} 0} &0\\ \hline
0111000& 0 & 0 & {\color{red} 0} & {\color{red} 0} & {\color{green} 0} & {\color{green} 0} &{\color{blue} 0} &{\color{blue} 1-1}& 0 \\ \hline
0111001& 0 & 0 & {\color{red} 0} & {\color{red} 0} & {\color{green} 0} & {\color{green} 0} &{\color{blue} 0} &{\color{blue} 1-1} & 0\\ \hline
0111010& 0 & 0 & {\color{red} 1-1} & {\color{red} 0} & {\color{green} 1-1} & {\color{green} 0} &{\color{blue} 0} &{\color{blue} 0} & 0\\ \hline
0111011& 0 & 0 & {\color{red} 1-1} & {\color{red} 1-1} & {\color{green} 0} & {\color{green} 0} &{\color{blue} 0} &{\color{blue} 0} & 0\\ \hline
0111100& 0 & 0 & {\color{red} 0} & {\color{red} 0} & {\color{green} 0} & {\color{green} 0} &{\color{blue} 0} &{\color{blue} 1-1} & 0 \\ \hline
0111101& 0 & 0 & {\color{red} 1-1} & {\color{red} 0} & {\color{green} 0} & {\color{green} 0} &{\color{blue} 0} &{\color{blue} 0} & 0\\ \hline
0111110& 0 & 0 & {\color{red} 0} & {\color{red} 0} & {\color{green} 0} & {\color{green} 0} &{\color{blue} 0} &{\color{blue} 0} & 0\\ \hline
0111111& 0 & 0 & {\color{red} 0} & {\color{red} 0} & {\color{green} 0} & {\color{green} 0} &{\color{blue} 0} &{\color{blue} 0} & 0\\ \hline
1001000& 1 & 1 & {\color{red} 0} & {\color{red} 0} & {\color{green} 0} & {\color{green} 0} &{\color{blue} 0} &{\color{blue} 0} & 0\\ \hline
1001001& 1 & 1 & {\color{red} 0} & {\color{red} 0} & {\color{green} 0} & {\color{green} 0} &{\color{blue} 0} &{\color{blue} 0} & 0\\ \hline
1001010& 0 & 1-1 & {\color{red} 0} & {\color{red} 0} & {\color{green} 1} & {\color{green} 1} &{\color{blue} 0} &{\color{blue} 0} & 0 \\ \hline
1001011& 0 & 1 & {\color{red} 0} & {\color{red} 1} & {\color{green} 0} & {\color{green} 1} &{\color{blue} -1} &{\color{blue} 0} & -1 \\ \hline
1001100& 1-1 & 1-1 & {\color{red} 0} & {\color{red} 0} & {\color{green} 0} & {\color{green} 0} &{\color{blue} 1} &{\color{blue} 1} & 0 \\ \hline
1001101& 0 & -1 & {\color{red} 1} & {\color{red} 0} & {\color{green} 0} & {\color{green} -1} &{\color{blue} 1} &{\color{blue} 0} &+1\\ \hline
1001110& 0 & 0 & {\color{red} 0} & {\color{red} 0} & {\color{green} 0} & {\color{green} 0} &{\color{blue} 1-1} &{\color{blue} 0} & 0  \\ \hline
1001111& 0 & 0 & {\color{red} 0} & {\color{red} 0} & {\color{green} 0} & {\color{green} 0} &{\color{blue} 1-1} &{\color{blue} 0} & 0\\ \hline
1011000& -1 & 0 & {\color{red} 0} & {\color{red} 1} & {\color{green} -1} & {\color{green} 0} &{\color{blue} 0} &{\color{blue} 1} & -1\\ \hline
1011001& -1 & 0 & {\color{red} 0} & {\color{red} 1} & {\color{green} -1} & {\color{green} 0} &{\color{blue} 0} &{\color{blue} 1} & -1\\ \hline
1011010& 0 & 0 & {\color{red} 1} & {\color{red} 1} & {\color{green} 1-2} & {\color{green} -1} &{\color{blue} 0} &{\color{blue} 0} & 0 \\ \hline
1011011& 0 & 0 & {\color{red} 1} & {\color{red} 2-1} & {\color{green} -1} & {\color{green} -1} &{\color{blue} 0} &{\color{blue} 0} & 0\\ \hline
1011100& 0 & 0 & {\color{red} 0} & {\color{red} 1-1} & {\color{green} 0} & {\color{green} 0} &{\color{blue} 0} &{\color{blue} 1-1} & 0\\ \hline
1011101& 0 & 0 & {\color{red} 1-1} & {\color{red} 1-1} & {\color{green} 0} & {\color{green} 0} &{\color{blue} 0} &{\color{blue} 0} & 0\\ \hline
1011110& 0 & 0 & {\color{red} 0} & {\color{red} 1-1} & {\color{green} 0} & {\color{green} 0} &{\color{blue} 0} &{\color{blue} 0} & 0 \\ \hline
1011111& 0 & 0 & {\color{red} 0} & {\color{red} 1-1} & {\color{green} 0} & {\color{green} 0} &{\color{blue} 0} &{\color{blue} 0} & 0\\ \hline
1101000& 1 & 0 & {\color{red} 1} & {\color{red} 0} & {\color{green} 1} & {\color{green} 0} &{\color{blue} 0} &{\color{blue} -1} & +1 \\ \hline
1101001& 1 & 0 & {\color{red} 1} & {\color{red} 0} & {\color{green} 1} & {\color{green} 0} &{\color{blue} 0} &{\color{blue} -1} & +1 \\ \hline
1101010& 0 & -1 & {\color{red} 1} & {\color{red} 0} & {\color{green} 2} & {\color{green} 1} &{\color{blue} 0} &{\color{blue} -1} & +1 \\ \hline
1101011& 0 & 0 & {\color{red} 1} & {\color{red} 1} & {\color{green} 1} & {\color{green} 1} &{\color{blue} -1} &{\color{blue} -1} & 0 \\ \hline
1101100& -1 & 0 & {\color{red} 1-1} & {\color{red} 1} & {\color{green} -1} & {\color{green} 0} &{\color{blue} 0} &{\color{blue} 1}&-1  \\ \hline
1101101& 0 & 0 & {\color{red} 2-1} & {\color{red} 1} & {\color{green} -1} & {\color{green} -1} &{\color{blue} 0} &{\color{blue} 0} & 0\\ \hline
1101110& 0 & 0 & {\color{red} 1-1} & {\color{red} 1-1} & {\color{green} 0} & {\color{green} 0} &{\color{blue} 0} &{\color{blue} 0} & 0\\ \hline
1101111& 0 & 0 & {\color{red} 1-1} & {\color{red} 1-1} & {\color{green} 0} & {\color{green} 0} &{\color{blue} 0} &{\color{blue} 0} & 0\\ \hline
1111000& 0 & 0 & {\color{red} 0} & {\color{red} 0} & {\color{green} 0} & {\color{green} 0} &{\color{blue} 0} &{\color{blue} 1-1} & 0\\ \hline
1111001& 0 & 0 & {\color{red} 0} & {\color{red} 0} & {\color{green} 0} & {\color{green} 0} &{\color{blue} 0} &{\color{blue} 1-1} & 0\\ \hline
1111010& 0 & 0 & {\color{red} 1-1} & {\color{red} 0} & {\color{green} 1-1} & {\color{green} 0} &{\color{blue} 0} &{\color{blue} 0} & 0\\ \hline
1111011& 0 & 0 & {\color{red} 1-1} & {\color{red} 1-1} & {\color{green} 0} & {\color{green} 0} &{\color{blue} 0} &{\color{blue} 0} & 0\\ \hline
1111100& 0 & 0 & {\color{red} 0} & {\color{red} 0} & {\color{green} 0} & {\color{green} 0} &{\color{blue} 0} &{\color{blue} 1-1} & 0 \\ \hline
1111101& 0 & 0 & {\color{red} 1-1} & {\color{red} 0} & {\color{green} 0} & {\color{green} 0} &{\color{blue} 0} &{\color{blue} 0} & 0 \\ \hline
1111110& 0 & 0 & {\color{red} 0} & {\color{red} 0} & {\color{green} 0} & {\color{green} 0} &{\color{blue} 0} &{\color{blue} 0} & 0\\ \hline
1111111& 0 & 0 & {\color{red} 0} & {\color{red} 0} & {\color{green} 0} & {\color{green} 0} &{\color{blue} 0} &{\color{blue} 0} & 0 \\ \hline

	\caption{$\Delta$Difference, $d_n=1$}
	\label{tab:ecaj1}
\end{longtable}

\section{Conclusion}

The theorem allows us to answer the question by illustrating that for each binary sequence $D$ the property is verified.
This theorem can be generalized in order to prove (or to refuse) the same properties for inverse sequences of length bigger than 4.





\end{document}